\documentclass[a4paper]{article}

\author{Joaquim Ro\'e}

\usepackage[centertags]{amsmath}
\usepackage{amscd}
\usepackage{amsthm}
\usepackage{amssymb}
\usepackage[dvips]{graphics}

\theoremstyle{definition}
\newtheorem{Def}{Definition}[section]
\theoremstyle{plain}
\newtheorem{Lem}[Def]{Lemma}
\newtheorem{Cor}[Def]{Corollary}
\newtheorem{Pro}[Def]{Proposition}
\newtheorem{Teo}[Def]{Theorem}

\theoremstyle{remark}
\newtheorem{Rem}[Def]{Remark}

\makeatletter
\@addtoreset{table}{section}
\makeatother

\newcommand{\Ker}{\operatorname{Ker}}
\newcommand{\ord}{\operatorname{ord}}

\newcommand{\Hilb}{\operatorname{Hilb}}

\newcommand{\len}{\operatorname{length}}
\newcommand{\red}{\operatorname{red}}

\newcommand{\Tr}[4]{{\operatorname{Tr}_{#3}^{#4}(#1|#2)}}
\newcommand{\Res}[4]{{\operatorname{Res}_{#3}^{#4}(#1|#2)}}

\newcommand{\tr}[4]{{\operatorname{tr}_{#3}^{#4}(#1|#2)}}
\newcommand{\res}[4]{{\operatorname{res}_{#3}^{#4}(#1|#2)}}

\newcommand{\Cl}{{\mathit{Cl}}}

\newcommand{\C}{{\bf C }}

\newcommand{\M}{\mathfrak m}
\renewcommand{\C}{\mathbb{C}}
\newcommand{\Z}{\mathbb{Z}}

\newcommand{\m}{{\bf m}}

\renewcommand{\P}{\mathbb{P}}

\renewcommand{\O}{{\cal O}}

\begin{document}

\title{Maximal rank for planar singularities of multiplicity 2}
\maketitle

\begin{abstract}
We prove that general unions of singularity schemes of multiplicity two in the
projective plane have maximal rank.
\end{abstract}

\section{Introduction}

It has been known for a long time that given a set of $r\ne 2,5$ points in
general position in the plane and a positive integer $d$, the linear
system of all curves of degree $d$ singular at the $r$ points has
dimension $\max \{-1,d(d+3)/2-3r\}$. In other words, each singularity
imposes 3 linearly independent conditions or there are no curves with
the required singularities, except in the aforementioned cases. Early
references to this result can be traced back to Palatini \cite{Pal03}. The
standard modern reference is Hirschowitz \cite{Hir85}. 

General curves in the linear systems just described have ordinary
nodes as their only singularities; it is natural to ask about linear
systems of curves with more complicated singular points, and under
which hypotheses the conditions imposed by the singularities are
linearly independent.
By general principles it is clear that they will be independent if the
degree is large enough, and there are recent results due to Shustin
and his collaborators that provide bounds for what ``large enough''
must mean if the ``position'' of the singularities is general. These
bounds of \cite{GLS98} and \cite{Shu04} are valid for any type of
singularity, but they are not sharp, and for some kinds of
singularities (such as nodes, as above) it is possible to do better.

D. Barkats proved in \cite{Bar98} that the linear system of all curves of a
given degree with $\nu$ ordinary nodes and $\kappa$ ordinary cusps at given
(general) points and with given (general) tangents for the cusps has
dimension $\max \{-1,d(d+3)/2-3\nu -5\kappa\}$, except in the two
cases already encountered (2 or 5 nodes) and when there are two cusps.
In other words, also in this case each singularity imposes linearly
independent conditions or there are no curves with the required
singularities. The same is true for node and tacnode singularities
\cite{Roe01a}, except when the orders of the nodes and tacnodes add up
to 2 or 5 (they are coalescences of the classical 2-node and 5-node
cases). In Arnold's notation, this means that a collection of 
singularities of types $A_1$ and $A_2$ imposes independent conditions,
as does a collection of singularities of types $A_k$ where every $k$ is
odd. In this work (theorem \ref{lastthm}) we prove that this holds in
fact for every collection 
of singularities of multiplicity two (i.e., of types $A_k$ with $k$
arbitrary) with the only exceptions already known. 

\medskip

In order to precisely state what ``imposing a singularity of a given
type in a given position'' and ``general position'' mean we need some
algebro-geometric language. Let us fix the setting first. Except if
otherwise stated, we work over an algebraically closed field $k$ of
arbitrary characteristic. A
type of singularity means an equivalence class of germs 
of plane curve under equisingularity --two singular points are
equisingular if their embedded resolutions have the same
combinatorics. If $k=\C$ is the 
complex field, then equisingularity is the same as topological
equivalence (in a neighbourhood of the singular point). We work in the
projective plane $\P^2=\P^2_k$ although, as we deal with points in
general position, it is more or less indifferent to use an affine or
projective setting. 

The embedded resolution of a singular point of a plane curve consists
in blowing up the point and all singular points of its successive
strict transforms, so that at the end of the (finite) process one gets a
surface in which the strict transform of the curve is nonsingular and its
total transform is a normal crossings divisor. To the \emph{cluster} of
points that have been blown up we associate a combinatoric invariant, the
weighted \emph{Enriques diagram}, which is a tree whose vertices
represent the points, whose edges represent their proximity relations
and comes with integral weights, that represent their multiplicities
--a point is proximate to another if it lies in (the strict transform of) its
exceptional divisor. Two singular points are equisingular if and only
if the Enriques diagrams of the associated clusters coincide. 
We usually denote clusters by capital letters as $K$, whereas the
weights are denoted by $\m$, an Enriques diagram is $D$ and a weighted
diagram is $(D,\m)$.

The reader may find the basics on clusters and Enriques diagrams
in \cite{Cas00}. Here we need mainly the \emph{complete ideals}
defined by weighted clusters; let us briefly recall some of the basic
facts concerning them. Let $O$ be a point in the plane where
the curve $C$ has a singularity. Let $\O$ be the (two-dimensional,
regular) local ring at $O$ (so the germ of $C$ is defined as $f=0$ for
some $f\in \O$), and let $K$ be the cluster of the
embedded resolution of $C$. Then the set $I_{(K,\m)}$ of all $g\in \O$
such that 
the germ of curve $g=0$ goes through the points of $K$ with (virtual)
multiplicities at least as big as those of $C$ is an
ideal, which is $\M$-primary ($\M$ being the maximal ideal of $\O$)
and complete. For general $g\in I_{(K,\m)}$, $g=0$ is equisingular to $C$.
The equisingularity type is determined by the class of $g$ modulo
$\M^n$ for some $n$, and $I_{(K,\m)}$ is $\M$-primary, so one may use
$\O$ or its completion with respect to the maximal ideal, if it
simplifies matters.

\begin{figure}
  \begin{center}
    \mbox{\includegraphics{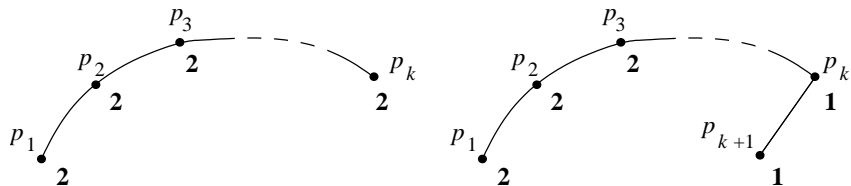}}
    \caption{Enriques diagrams of a tacnode (left, $A_{2k-1}$) and a
      cusp (right, $A_{2k-2}$). The multiplicities (weights) are in
      boldface. In both cases every point is
      proximate to its predecessor, but the last point $p_{k+1}$ of a
      cusp is proximate to $p_{k-1}$ as well, so it is a
      \emph{satellite}; this is represented in the diagram by a
      straight segment. These are all the singularities
      of multiplicity 2.}
    \label{fig:Ak}
  \end{center}
\end{figure}

If $C$ is a curve with several singular points, we associate to it a
cluster that is the disjoint union of the clusters of all its singularities,
and the corresponding Enriques diagram (which is now a forest rather
than a tree).
The set $\Cl(D)$ of all clusters with the same Enriques diagram $D$ has a
natural structure of quasiprojective algebraic variety \cite{Roe04b};
whenever we state some claim about singularities of type $(D,\m)$
in general position we mean that the claim
holds for singularities whose cluster lies in a Zariski open set of $\Cl(D)$.
Figure \ref{fig:Ak} shows the Enriques diagrams that appear for
singularities of multiplicity 2.

\medskip

Let $K$ be a cluster of points in $\P^2$ with some weights that
correspond to a singularity type.
Let $I_{(K,\m)}$ denote the ideal sheaf supported at the
proper points of $K$  which is locally defined as above, by the
condition of going through the points of $K$ with the assigned
multiplicities, and let $Z_{(K,\m)}$ denote the (zero-dimensional) subscheme
of $\P^2$ defined by $I_{(K,\m)}$. For every positive $d\in \Z$, the twisted
global sections $\Gamma(I_{(K,\m)}(d))$ are homogeneous polynomials of degree
$d$ defining curves that go through the points of $K$ with the
assigned multiplicities, and $\P(\Gamma(I_{(K,\m)}(d)))$ is the linear system
of curves of degree $d$ with the assigned singularities at the
assigned positions. If $d$ is high enough, then general curves in
$\P(\Gamma(I_{(K,\m)}(d)))$ do indeed have the singularity type given by the
Enriques diagram of $K$. The conditions imposed by $K$ are independent
or there are no curves of degree $d$ containing $Z_{(K,\m)}$ if and
only if the canonical map 
$$
k[x,y,z]_d\cong\Gamma(\O_{\P^2}(d)) \longrightarrow \Gamma(\O_{Z_{(K,\m)}}(d))
\cong \Gamma((\O_{\P^2}/I_{(K,\m)})(d))
$$
is either surjective or injective, i.e., has maximal rank. If this
happens for a given $(K,\m)$ and for all $d$ we say that $Z_{(K,\m)}$
has maximal 
rank. Thus for instance the first result mentioned in the introduction
may be rephrased in more fancy words by saying that, if $D$ is the
diagram consisting of $r$ unconnected vertices --no edges--, $K$ is general
in $\Cl(D)$, and we take all weights equal to 2 then $Z_{(K,\m)}$ has maximal
rank, and we claim that the same is true if $(D,\m)$ is a union of
weighted diagrams of the types shown in figure \ref{fig:Ak}.

It may be good to warn that the schemes $Z_{(K,\m)}$ obtained with a fixed
weighted Enriques diagram $(D,\m)$ need not be isomorphic in
general. The class of schemes isomorphic to a given $Z_{(K,\m)}$ is contained
in the class of all schemes $Z_{(K,\m)}$, where $K$ has diagram $(D,\m)$;
for almost all $(D,\m)$ there is a nontrivial moduli space of such
schemes. However, for types $(D,\m)$ of multiplicity two all $Z_{(K,\m)}$ are
indeed isomorphic (because $A_k$-singularities have no moduli, see
\cite{Arn76}) so in the sequel we seldom mention 
isomorphism classes of zero-dimensional schemes.

To show that general schemes of a given class have maximal rank it is
often useful to use \emph{specialization} and \emph{semicontinuity}:
if one has a flat family of schemes $Z_t$ parameterized by some smooth
scheme $t\in T$, such that for some special value of the parameter
$t=0$ the scheme $Z_0$ has maximal rank, then the principle of semicontinuity
\cite[chapter III, 12]{HAG} tells us that general members of the family have
maximal rank. The strategy of our proof consists in a sequence of
specializations which furnish a family $Z_t$ whose general members are
of type $Z_{(K,\m)}$ (where the Enriques diagram of K is a union of diagrams
of the types shown in figure \ref{fig:Ak}) and where $Z_0$ is known to
have maximal rank.

\medskip

Let $(D,\m)$ be a given union of
weighted diagrams of the types shown in figure \ref{fig:Ak}.
The first specialization simplifies matters by reducing to a family of
schemes supported at a single point. For every cluster $K$ with
diagram $D$ there is a smooth curve $C$ going through all the free
points of $K$ (i.e., through the subcluster $K'$ consisting
of every point of the tacnodes and every point of the cusps except the
last one, weighted with multiplicity 1 at all the points). It is
enough to pick $C$ of high enough degree and consider the 
complete (curvilinear) ideal associated to $K'$ with these weights. 
We allow the base points to move on $C$, and specialize them to ``collide''
at a single point, giving as a flat limit a zero-dimensional scheme
supported at a single point and contained in $2C$, the one-dimensional
scheme whose equation is the square of that of $C$. 

The specialization used leads to schemes that are not singularity
schemes, because they are not defined by complete ideals; their
defining property is to be contained in the double of a smooth
curve, and we call such schemes 2-curvilinear schemes (see section
\ref{sec:hilbscheme}). Actually we
prove a maximal rank statement for 2-curvilinear schemes; as far as
we know, this is the first place where a maximal rank result is proved
that involves schemes whose defining ideals are not complete. This can
be understood as a generalization of the well-known fact that general
curvilinear schemes have maximal rank (see \cite{CM98} in arbitrary
dimension, or apply \cite{Bri77} in dimension 2).

There are two numerical invariants naturally associated to a
2-curvilinear scheme $Z$. The first is the length $N$ of $Z$, and the
second is the maximal contact $\ell$ of $Z$ with a smooth
curve whose double contains $Z$; they satisfy the inequalities 
$0\le \ell \le N \le 2\ell$. Our main result is the following.

\begin{Teo}
\label{mainthm}
Let $N, \ell$ be two positive integers with $0\le \ell \le N \le 2\ell-1-3
\sqrt{N-\ell}$. Then for every isomorphism class of 2-curvilinear
schemes $Z$ whose length is $N$ and whose maximal contact with smooth
curves whose double contains $Z$ is $\ell$, general members of the
class in $\P^2$ have maximal rank.
\end{Teo}

This statement makes sense because zero-dimensional
schemes of given isomorphism class form an irreducible family.
In section \ref{sec:2curv} we prove the theorem, and we also give 
a precise description of the kind of schemes obtained as
limits in the collision above. Once this is done, we get as a
corollary the following.

\begin{Teo}
\label{secondthm}
Let $(D,\m)$ be a union of singularity types of multiplicity two
(whose weighted diagrams are of the types shown in
figure \ref{fig:Ak}). Then singularity schemes of type $(D,\m)$ in
general position have maximal rank in degrees $d\ge 13$.
\end{Teo}

To give a complete proof of the result, theorem \ref{lastthm},
claimed above we need
to deal with the (finitely many) cases that involve degrees 12 or
less, which we do case by case with ad-hoc methods.

Most of the paper is devoted to the proof of theorem \ref{mainthm},
which is done --as said-- by providing a sequence of
specializations. The specializations are relatively easy to describe,
and the main difficulty relies in computing the limit of an explicit
one-parameter family of zero-dimensional schemes. To do this we rely
in an algebraic lemma in the spirit of Alexander and Hirschowitz
\cite{AH00} or \'Evain \cite{Evain}.

\section{An algebraic lemma}
\label{sec:formal}

In this section, we state and prove a slightly generalized version of
the ``differential Horace lemma'' of \cite[section 8]{AH00}. 
The generalization, which is quite natural and more or less implicit
in the works of Alexander-Hirschowitz \cite{AH92, AH95, AH00}, \'Evain
\cite{Eva98, Eva99}, Mignon
\cite{Mig00a,Mig01} and others,
allows us to deal with more general families of zero-dimensional schemes, whose
defining ideals are not ``vertically graded'', as required in
the cited papers.

Let $R$ be an integral $k$-algebra, and consider $R_t=R \otimes
k[[t]]$. Given $f_t\in R_t$, denote $f_0\in R$ its image by the
morphism $t \mapsto 0$. Similarly, for an ideal $I_t$ in
$R_t$, denote $I_0=(I_t+(t))/(t) \subset R_t/(t) \cong R$.

Given an ideal $I_t \in R_t$, an element $y \in R$ and an integer
$p\ge 1$, the $p$-trace and $p$-residual ideals of $I_t$
with respect to $y$ are defined as follows:
\begin{align*}
\Tr{I_t}{y}{p}{}= & \frac{\left((I_t+(y)):t^{p-1}\right)_0}{(y)}
\subset R/(y),\\ \Res{I_t}{y}{p}{}= & \left((I_t+(t^{p})):y\right)_0
\subset R.
\end{align*}
Note that there are inclusions $\Tr{I_t}{y}{1}{}\subset \Tr{I_t}{y}{2}{}
\subset \dots$, and $\Res{I_t}{y}{1}{}\supset \Res{I_t}{y}{2}{} \supset
\dots I_0$. 
The ideals we are interested in have generically finite colength;  we
define $\tr{I_t}{y}{p}{}=\dim_k((R/(y))/\Tr{I_t}{y}{p}{})$ and
$\res{I_t}{y}{p}{}=\linebreak[4]
\dim_k(R/\Res{I_t}{y}{p}{})$.

As in \cite{AH00}, given any linear subspace $V \subset R$ and $y\in R$, let
$\Res{V}{y}{}{}=\{v \in R \ | \ vy \in V\}$. Since $y$ is not a zero-divisor, we get a
residual exact sequence
$$ 0 \longrightarrow \Res{V}{y}{}{} \overset{y}\longrightarrow V
\longrightarrow V/V \cap (y) \longrightarrow 0. $$
\begin{Pro}
\label{formal} Let $V \subset R$ be a $k$-linear subspace, and $I_t
\subset R_t$ an ideal such that $R_t/I_t$ is flat over 
  $k[[t]]$. Let $p\in \Z$ and $y \in R$ be given, with $p \ge 1$.
Consider the following three canonical maps:
$$ 
\frac{V}{V \cap (y)} \overset{\varphi_p}{\longrightarrow}
\frac{R/(y)}{\Tr{I_t}{y}{p}{}}, \quad
\Res{V}{y}{}{} \overset{{\check\varphi}_p}{\longrightarrow}
\frac{R}{\Res{I_t}{y}{p}{}}, \quad
V \otimes k[[t]] \overset{\varphi_t}{\longrightarrow} 
R_t/I_t. $$
If $\varphi_p$ is injective, then $(\Ker \varphi_t)_0 \subset y \,
\Ker {\check\varphi}_p$.
\end{Pro}

\begin{proof}
Let $f_t \in \Ker \varphi_t= V \otimes k[[t]] \cap I_t$. 
If $f_t\in (t^{p}, y)$,
i.e.,  $f_t=g_t y+h_t t^p$ for some $g_t,h_t\in R_t$ 
then by the definitions
$g_0 \in \Res{V}{y}{}{}\cap \Res{I_t}{y}{p}{}=
\Ker {\check\varphi}_p$, 
and therefore $f_0=y g_0 \in y \Ker {\check\varphi}_p$, so it will be
enough to prove that the injectivity of  $\varphi_p$
implies $f_t\in (t^{p}, y)$.

Write $f_t=\sum F_jt^j$, with $F_j \in R$. Denote $\bar F_j$ the class
of $F_j$ in $R/(y)$; we 
want to see that $\bar F_0= \dots= \bar F_{p-1}=0$.

The inclusions $\Tr{I_t}{y}{0}{}\subset
\Tr{I_t}{y}{1}{}  \subset \dots$
together with the injectivity of  $\varphi_p$ tell us that,
for every $j=1, \dots, p$, the map $$
\varphi_j:\frac{V}{V\cap (y)}
\longrightarrow
\frac{R/(y)}{\Tr{I_t}{y}{j}{}}$$ is
injective. As we have $f_t \in I_t$, it follows
that $$ \bar F_0 \in \frac{I_t+(y,t)}
{(y,t)}=\Tr{I_t}{y}{1}{}, $$ i.e., $\varphi_1 (\bar
F_0)=0$, and therefore $\bar F_0=0$. Now we argue by iteration:
let $1\le i< p$, and assume we have proved $\bar F_0= \dots= \bar
F_{i-1}=0$. This means that $f_t \in (y, t^i)$, so
$\sum_{j\ge i}F_j t^{j-i} \in
(I_t+(y)):t^i$, which implies $\bar F_i \in
\Tr{I_t}{y}{i+1}{},$ i.e., $\varphi_{i+1}
(\bar F_i)=0$, and therefore $\bar F_i=0$. The proof is now
complete.
\end{proof}

\section{Monomial ideals}
\label{sec:monomial}

In this paper, the ring $R$ above will be the completion of the local
ring at a given point of a smooth algebraic surface, and therefore isomorphic
to a power series ring $R\cong k[[x, y]]$. In particular it is a
regular local ring, and $R_t\cong k[[x,y,t]]$ is a regular local ring
as well. Their maximal ideals are 
$\mathfrak{m}=(x, y)$ and $\mathfrak{m}_t=(x,y,t)$ respectively.

We shall be dealing with a restricted kind of ideals, 
of the form
$$
I_E=\left(x^{e_1}f^{e_2}\right)_
{(e_1,e_2)\in E},
$$
where $E \subset \Z_{\ge 0}^2$ is a \emph{staircase}, that is, 
$E + \Z_{\ge 0}^2\subset E$ (see \cite{Eva??}), and
$f=x+y+t$. More generally, in some instances we shall consider ideals
of the form  
$$
I_{(E,f,g)}=\left(g^{e_1}f^{e_2}\right)_
{(e_1,e_2)\in E},
$$
where $f,g \in R$ are arbitrary.
For convenience, we introduce some language to deal with the
combinatorics of staircases.

\begin{Def}
  If $E\subset \Z_{\ge 0}^2$ is a staircase, we say that the
  \emph{length} of its $i$th stair is $\ell_E(i)=\inf \{e\ |\ (e,i)\in
  E\}$, and the \emph{height} of its $i$th ``slice'' is $h_E(i)=\inf
  \{e\ |\ (i,e)\in E\}$. We shall use the first differences of $\ell$
  and $h$ as well: $\hat \ell_E(i)=\ell_E(i)-\ell_E(i+1)$, $\hat
  h_E(i)=h_E(i)-h_E(i+1)$ . 
\end{Def}

When $E$ is a staircase with finite complement, $h_E$, $\ell_E$, $\hat
\ell_E$ and $\hat h_E$ are functions $Z_{\ge 0} \rightarrow Z_{\ge
  0}$, and each of them determines $E$ uniquely.

\begin{Lem}
\label{finideterm}
  Let $E$ be a staircase, and let $\alpha, \beta \in \Z$ be such that
  $\hat \ell(i)\le \alpha$ for all $i$ with $\ell(i)\ne 0$ and
  $\hat h(i)\le \beta$ for all $i$ with $h(i)\ne 0$. Let $f, g, f',
  g' \in R$ be such that $(f,g)=(f',g')$, $f-f'\in (f,g)^\alpha$ and
  $g-g'\in (f,g)^\beta$. Then $I_{(E,f,g)}=I_{(E,f',g')}$.
\end{Lem}

\begin{proof}
  If $E$ is empty then $I_{(E,f,g)}=(0)=I_{(E,f',g')}$ and there is
  nothing to prove; otherwise the hypotheses on the lengths of its
  stairs imply that $E$ has finite complement.

  We are going to prove that if $\hat h(i)\le \beta$ for all $i$
  with $h(i)\ne 0$ and $g-g'\in (f,g)^\beta$ then
  $I_{(E,f,g')}\subset I_{(E,f,g)}$. Then it will follow that $g-g'\in
  (f,g')^\beta=(f,g)^\beta$ and hence $I_{(E,f,g)}\subset
  I_{(E,f,g')}$; thus $I_{(E,g,g')}= I_{(E,f,g)}$. By symmetry then
  it will follow that if $\hat \ell(i)\le \alpha$ for all $i$
  with $\ell(i)\ne 0$ and $f-f'\in (f,g)^\alpha$ then
  $I_{(E,f',g)}= I_{(E,f,g)}$. Finally, $I_{(E,f,g)}=I_{(E,f',g')}$ as claimed.

  Let us see that for all $(e_1, e_2) \in E$, $g'^{e_1}f^{e_2} \in
  I_{(E,f,g)}$. We do it 
  by induction on $e_1$. If $e_1=0$ then there is nothing to prove;
  assume that $e_1>0$ and $g'^{e'_1}f^{e'_2} \in  I_{(E,f,g)}$ for all 
  $(e'_1, e'_2) \in E$ with $e'_1<e_1$.

  Put $h=g'-g \in (f,g)^\beta$. Then
$$
g'^{e_1}f^{e_2}=\sum_{i=0}^{e_1} \binom{e_1}{i} g^i h^{e_1-i} f^{e_2}. 
$$
  As $h^{e_1-i}\in (f,g)^{\beta(e_1-i)}$, we shall be done if
  $g^{i+a}f^{e_2+b}\in I_{(E,f,g)}$ for every $0\le i \le e_1$ and
  every $a,b\ge 0$ with $a+b=\beta(e_1-i)$ or, equivalently, if 
  $(i+a,e_2+b)\in E$ for every $0\le i \le e_1$ and
  every $a,b\ge 0$ with $a+b=\beta(e_1-i)$.
  But it follows from the hypothesis on $\beta$ that $(a',b')\in E$
  whenever $b'\ge e_2$ and $\beta a'+b' \ge \beta e_1+e_2$, which is
  easy to check for $(a',b')=(i+a,e_2+b)$.
\end{proof}

The computation of quotient
ideals of monomial ideals $I_E$ as above leads, under suitable
conditions, to new monomial ideals obtained by slicing off part of the
staircase. This fact has already been exploited by
Alexander-Hirschowitz and \'Evain, and we shall take advantadge of it
as well. So define $\sigma(E,p)$ as the staircase obtained from $E$ by
deleting the $p$th slice, i.e., the unique staircase whose height
function has
\begin{equation*}
  h_{\sigma(E,p)}(i)=
  \begin{cases}
    h_E(i) & \text{if }i\le p, \\
    h_E(i+1) & \text{if }i > p.
  \end{cases}
\end{equation*}

\begin{Pro} 
\label{flatstairs}
  Let $E \subset \Z_{\ge 0}^2$ be a staircase,
  and $I_t=I_E$ as defined above. Assume that $\hat
  \ell_E(i)\ge 2$ for all $i<h_E(0)-1$. Then 
  \begin{enumerate}
  \item $R_t/I_t$ is flat over $k[[t]]$ and over $k[[y]]$,
  \item $\tr{I_E}{y}{p}{}=h_E(p-1)$,
  \item \label{restairs} $\Res{I_t}{y}{p}{}=(I_{\sigma(E,p)}+(t))/(t)$.
\end{enumerate}
\end{Pro}

\begin{proof}
  To prove the first claim, consider the automorphism $\psi$ of
  $R_t$ defined by 
  $\psi(x)=x$, $\psi(y)=f$, $\psi(t)=t$. It is a
  $k[[t]]$-automorphism (it leaves $k[[t]]$ fixed) so $R_t/I_t$ is
  flat over $k[[t]]$ if and only if $R_t/\psi^{-1}(I_t)$ is. But
  $\psi^{-1}(I_t)$ is generated by monomials in $x$ and $y$, so
  $R_t/\psi^{-1}(I_t)=R/\psi^{-1}(I_t)_0 \otimes k[[t]]$ is obviously
  flat over $k[[t]]$. The same argument, reversing the roles of $y$
  and $t$, proves the flatness over $k[[y]]$.

It has already been remarked and used (see \cite[section 8.1]{AH00},
\cite{Eva??}) and it is not hard to prove directly that for 
$$J_E=\left((y+t)^{e_1}x^{e_2}\right)_{(e_1,e_2)\in E},$$  
$\tr{J_E+(t^q)}{y}{p}{}=\tr{J_E}{y}{p}{}=h_E(p-1)$,
and $\Res{J_t}{y}{p}{}=(J_{\sigma(E,p)}+(t))/(t)$.

Consider now the automorphism $\psi$ of $R_t$ defined by 
$\psi(x)=f$, $\psi(y)=y$, $\psi(t)=t$. Leaving $t$ fixed, it
induces an automorphism of $R\cong R_t/(t)$ which we also denote by $\psi$.
It is not hard to see that the hypothesis $\hat \ell_E(i)\ge 1$ for
all $i<h_E(0)-1$ implies that $\hat h_E(i)\le 1$ for all $i$, and
therefore the previous lemma tells us that
$I_E=\left((y+t)^{e_1}f^{e_2}\right)_{(e_1,e_2)\in E}$, so
$\psi^{-1}(I_E)=J_E$, and therefore again
$\tr{I_E}{y}{p}{}=\tr{J_E}{y}{p}{}=h_E(p-1)$, as desired.

Similarly, we have $\Res{I_t}{y}{p}{}=
\psi(\Res{J_t}{y}{p}{})=\psi(J_{\sigma(E,p)})+(t)/(t)$.
As $\hat \ell_E(i)\ge 2$ for all $i<h_E(0)-1$, it follows immediately that $\hat \ell_{\sigma(E,p)}(i)\ge 1$ for
all $i<h_E(0)-1$, and therefore $\psi(J_{\sigma(E,p)})=
\left((y+t)^{e_1}f^{e_2}\right)_{(e_1,e_2)\in
  \sigma(E,p)}=I_{\sigma(E,p)}$, finishing the proof.
 \end{proof}

\section{Adjacencies in the Hilbert scheme}
\label{sec:hilbscheme}

Consider now a point $O\in \P^2$, its blowing-up 
$\pi_O:S \longrightarrow \P^2$, and a point $O'$ in the first
(infinitesimal) neighbourhood of $O$
(i.e., $O'\in S$ and $\pi_O(O')=O$). Let $\O_{O',S}$ be the local
ring at $O'$ and $y \in \O_{O',S}$ a local equation of the exceptional
divisor $D=\pi_O^{-1}(O)$. To every $f,g \in \O_{O',S}$
and every staircase $E$ we associate the ideal 
$I_{(E,f,g)}$ as defined in the previous section 
and for every integer $m$ its ($m$-twisted) push-forward
$$J_{(E,f,g,m)}=(\pi_O)_* \left(y^mI_{(E,f,g)}\right)\subset
\O_{O,\P^2}.$$
If $f$ and $g$ have no common components, and $E$ has
finite complement, then $J_{(E,f,g,m)}$ is $\M$-primary (where
$\M$ is the maximal ideal of $\O_{O,\P^2}$) and hence it defines a
zero-dimensional subscheme $Z_{(E,O,O',f,g,m)}$ of $\P^2$ 
supported at $O$.

We define \emph{2-curvilinear schemes} to be those schemes locally contained
in the double of a curve; more precisely, a zero-dimensional scheme
$Z$ is called 2-curvilinear if it satisfies the
following properties: 
\begin{enumerate}
\item $Z$ has embedding dimension at most 2; i.e., for every maximal ideal
  $\M$ of the Artinian ring $\O_Z$, $\dim_k \M/\M^2\le 2$, and
\item for every maximal ideal $\M$ of $\O_Z$, there
  exists $f\in \M \setminus \M^2$ with $f^2=0$. Such an $f$ is not
  unique, and we shall assume that it has been chosen of
  maximal contact with $Z$, i.e., that for every $g\in \M \setminus
  \M^2$ with $g^2=0$, $\dim_k (\O_Z/(f))_\M \ge \dim_k
  (\O_Z/(g))_\M$. 
\end{enumerate} 
Moreover, to every such $Z$ we attach invariants $N=\dim_k \O_Z$
(length) and $\ell$ (contact), which if $Z$ is irreducible can be
computed as $\ell=\dim_k(\O_Z/(f))$, where $f$ is
the chosen $f\in \M \setminus \M^2$ with $f^2=0$ (for the unique
maximal ideal). If $Z$ has several components then its invariants are
simply the sum of the invariants of each component.

Our interest in schemes $Z_{(E,O,O',f,g,m)}$ arises from the fact that
some of them 
are specializations of (unions of) singularity schemes of multiplicity
two. More precisely, they sit (inside $\Hilb \P^2$) in the closure of
the (irreducible) subscheme parameterizing 2-curvilinear schemes.
In order to see this, we begin by showing
that the $Z_{(E,O,O',f,g,m)}$ form
a ``nice'' subset of $\Hilb \P^2$. 
To simplify,  assume that $g=x\in \O_{O',S}$ is transverse to $D$
(i.e., $x,y$ are a system of parameters of $\O_{O',S}$), and denote 
$s=\ord (f|_D)=\dim_k((\O_{O',S}/(y,f))$ the
intersection multiplicity of $f=0$ with the exceptional divisor. 
Introduce the notation
$$
H_{m,E,s}=\left\{Z_{(E,O,O',f,g,m)}\left| 
  \begin{matrix}
    O \in \P^2; \, \pi_O:S_O \rightarrow \P^2 \text{ is the blowing-up
    of }O; \hfill \\
    O'\in D=\pi_O^{-1}(O); \,  f,g \in \O_{O',S_O}; \, \hfill\\
    (f,g)=\M_{O',S_O};\, \ord (g|_D)=1;\, \ord (f|_D)=s. \hfill
  \end{matrix}
\right. \right\}
$$

 \begin{Lem}
 \label{codimform}
    Let $O\in \P^2$, $O'\in S$, $x,y,g\in\O_{O',S}$ be
   given as above. Put $s=\ord (f|_D)$, let $E$ be a staircase
   with finite complement, and let $m_0=\min \{e_1+se_2|(e_1,e_2)\in
   E\}$. Then for every integer $m\ge m_0-1$, 
 $$\dim_k \frac{\O_{\P^2,O}}{J_{(E,x,f,m)}}=
 \len(m,E):=\binom{m+1}{2} + \#(\Z_{\ge 0} \setminus E).$$
 \end{Lem}
 \begin{proof}
   Let $X\subset S$ be the zero-dimensional scheme defined by
   $I_{E,f,x}$ (so $X_{\red}=O'$). It is clear that $\len X=\#(\Z_{\ge
   0} \setminus E)$ and, denoting by $D=\pi_O^{-1}(O)$ the exceptional
   divisor, $\len (X\cap D)=m_0$. Thus
   the claim follows by \cite[2.14]{Del73}. It is also possible to
   prove it along the lines of \cite[4.7.1]{Cas00}.
 \end{proof}

\begin{Lem}
Let $E$ be a staircase of height two and finite complement, and let
  $m_0=\min \{e_1+se_2|(e_1,e_2)\in E\}$. 
For every integer $s \ge 1$ and $m\ge m_0-1$, the set
$$
H_{m,E,s}\subset \Hilb^{\len(m,E)} \P^2,
$$
is an irreducible constructible subset for the Zariski topology.
\end{Lem}
\begin{proof}
Let $N=\len(m,E)$. The claim will follow from the existence of a
morphism $X \rightarrow \Hilb^{N} \P^2$ where $X$ is a smooth quasiprojective
variety, whose image is $H_{m,E,s}$. 

Assume that $\hat\ell_E(0)>0$. Recall from \cite{Roe04b} that the set
of all clusters with given 
Enriques diagram has a natural structure of quasiprojective variety. 
Let $r=1+\max(\hat\ell_E(0),\hat\ell_E(1))$, and
consider the diagram $D_s(r)$ shown in figure \ref{fig:proxo}. 
Take $X=\Cl(D_s(r))$
and define the map $Z:X \rightarrow \Hilb^{N} \P^2$ as follows.

\begin{figure}
  \begin{center}
    \mbox{\includegraphics{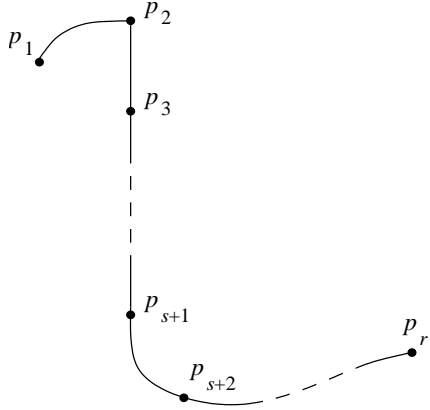}}
    \caption{Every point is proximate to the previous one; in
      addition, the $s$ points after the origin are proximate to it
      (so $p_3, \dots, p_{s+1}$ are satellite).}
    \label{fig:proxo}
  \end{center}
\end{figure}

Given $K\in X$, let $O=p_0(K)$, $O'=p_1(K)$, and remark that $O'$ is
in the first neighbourhood of $O$. Let $y\in \O_{O',S_O}$ be a local
equation of the exceptional divisor of blowing up $O$, and choose a
transverse germ $x=0$ not going through $p_2(K)$. Choose $f\in
\O_{O',S_O}$ to be
a local equation of a germ of curve smooth at $O'$ and going through
all points of $K$. Then we set $Z(K)=Z_{(E,O,O',f,x,m)}$. Note that by
the assumption that $\hat \ell_E(0)>0$ we have $\hat h_E(i)\le 1$ for
all $i$ and therefore by lemma \ref{finideterm} $I_{(E,f,x)}$ does not
depend on the choice of $x$, so neither does $Z(K)$.
Similarly, the definition of $r$ guarantees that $Z(K)$ does not
depend on the choice of $f$. 

It remains to be seen that the constructed map $Z:X \rightarrow
\Hilb^{N} \P^2$  is algebraic, and it is enough to do it locally. 

Let $K_0\in X$ be a closed point, and $p_1(K_0)\in \P^2$ the base point of
the corresponding cluster. If $(u,v)$ are affine coordinates in a
neighbourhood $U_0$ of $p_1(K_0)\in \P^2$, we may choose coordinates
$(u,v,x,y)$ in a neighbourhood $U_1$ of $p_2(K_0)$ in the variety $X_1$ of all
clusters of two points, in such a way that
\begin{align*}
U_1&\overset{\psi_1}{\rightarrow} U_0 & U_1&\overset{\pi_1}{\rightarrow} U_0\\
(u,v,x,y)&\mapsto (u,v)& (u,v,x,y)&\mapsto (u+x,v+xy)
\end{align*}
are local expressions of the structure morphism and the relative
blowing-up morphism of \cite{Roe04b}, i.e., $\pi_1$
restricted to the fiber of $\psi_1$ over $p\in U_0 \subset \P^2$ is
(an affine chart of) the blowing up of $p$, and $y=0$ is a local
equation of the relative exceptional divisor, i.e., its restriction to
each fiber of $\psi_1$ is a local equation of the corresponding
exceptional divisor.

With these coordinates, there is a function 
$f=y+a_sx^s+\dots + a_{r-1}x^{r-1} \in \O(X_1)$ whose restriction to
the fiber of $\psi_1$ over $p_1(K_0)$ vanishes
at all points of $K_0$ \cite{Cas00}, and therefore  
$Z(K_0)=Z_{(E,p_1(K_0),p_2(K_0),f,x,m)}$.

Then there are affine coordinates $(u,v,\bar x,x_s,\dots,x_{r-1})$ in a
neighbourhood $U_{r-1}$ of $K_0$ in $\Cl(D_s(r))$ such that the
restriction of 
$$
\tilde f=y+(a_s+x_s)(\bar x-x)^s+ \dots + (a_{r-1}+x_{r-1})(\bar x
-x)^{r-1} \in \O(U_{r-1}\times_{U_0} U_1)
$$
to the fiber over $p_1(K)$ is a local equation of a curve
going through all the points of $K$ (if $K\in U_{r-1}$ has
coordinates $(u,v,\bar x,x_s,\dots,x_{r-1})$). 

Consider now the ideals
\begin{align*}
 \mathcal{I}_{(E,\tilde f,x)}&=\left(x^{e_1}\tilde f^{e_2}\right)_
{(e_1,e_2)\in E} \subset \O(U_{r-1}\times_{U_0} U_1)\\
 \mathcal{J}_{(E,\tilde f,x,m)}&=(\text{id}\times\pi_1)_* 
\left(y^m\,\mathcal{I}_{(E,\tilde f,x)}\right)\subset
\O({U_{r-1}}\times U_0).
\end{align*}
Observe that (as in the previous lemma) for every $K\in U_{r-1}$,
  $m_0$ is the length of the 
  intersection (in the corresponding fiber) of the zeroscheme defined
  by the restriction of $\mathcal{I}_{(E,\tilde f,x)}$ with the
  exceptional divisor, so 
  by \cite[2.14]{Del73} and \cite[7]{EGA3}, $\O({U_{r-1}}\times U_0) /
\mathcal{J}_{(E,f,x,m)}$ is flat over $\O({U_{r-1}})$ of relative
length $N$, and it determines a morphism $U_{r-1}\rightarrow
\Hilb^N{U_0}\subset \Hilb^N{\P^2}$ which is set-theoretically equal to
  $Z$ above. 

It remains to deal with the case $\hat\ell_E(0)=0$, in which 
the ideal $I_{(E,f,x)}$ does depend on the choice of $x$ (in fact, on
the class of $x$ modulo the square of the maximal ideal of
$p_2(K)$). This choice can be 
parameterized by the set of free points in the first neighbourhood of
$p_2$. Thus one gets a map $X \rightarrow \Hilb^N \P^2$ as before,
with $X=\Cl(D_s(r))\times_{X_1}\Cl(D_1(3))$, which can be shown to be
a morphism in the same way. We leave the details to the interested reader.
\end{proof}

\begin{Teo}
\label{espbloc}
  Let $E$ be a staircase of height two and $s$ a positive integer
  satisfying $\hat \ell_E(0)\ge s+2$ and $\ell_E(1) \ge s$. Define
  $E_1$ to be the unique staircase of height (at most) two with
 $\ell_{E_1}(0)= \ell_E(0)-s-1$, $ \ell_{E_1}(1)=\ell_E(1)-s$. If
  $\ell_E(1) \ge 2s+1$, define furthermore $E_2$ to be the unique
  staircase of height (at most) two with 
 $\ell_{E_2}(0)= \ell_E(0)-2s-2$, $ \ell_{E_1}(1)=\ell_E(1)-2s-1$.

If $\ell_E(1) \le 2s$ then $H_{2s+1,E_1,s+1}\subset
\overline{H_{2s,E,s}},$ and if $\ell_E(1) \ge 2s+1$ then 
$H_{2s+2,E_2,s+1}\subset \overline{H_{2s,E,s}}.$
\end{Teo}

\begin{proof} Let $i=1$ if $\ell_E(1) \le 2s$ and $i=2$ if
  $\ell_E(1) \ge 2s+1$. Let $Z\in H_{m+i,E_i,s+1}$ be given by the ideal 
$$J=J_{(E_i,f,x,m+i)}=(\pi_O)_* \left(y^{m+i}I_{(E_i,f,x)}\right)\subset
\O_{O,\P^2},$$
where $O$ is some point in $\P^2$, $O'$ is a point in the first
neighbourhood of $O$, $x,y,f \in \O_{O',S_O}$ are smooth germs, $x,y$
are local parameters, $y=0$ is a 
local equation of the exceptional divisor $D$, and $\ord(f|_D)=s+1$.
Consider $f_t=f+tx^s\in \O_{O',S_O}\otimes k[t]$. For values of $t$ in
a neighbourhood of $0$, $J_t=(\pi_O)_* \left(y^m
  I_{(E,f,x)}\right)\subset \O_{O,\P^2}$ defines a zero-dimensional scheme $Z_t$ in
$H_{2s,E,s}$, so if we see that $Z$ is the flat limit of $Z_t,
t\rightarrow 0$, we shall be done. By \ref{codimform}, $Z$ and $Z_t$
have the same length, so it will be enough to show that $Z \supset
\lim Z_t$ or, equivalently, that for every $g_t \in J_t$, $g_0 \in J$.

Denote by $W\subset \O_{O',S_O}$ the set of virtual transforms of
equations of germs with (virtual multiplicity) at least $m$ at $O$,
i.e., $W=\pi^*(\M^m)/y^m$. If we show that for every $g_t \in I_t \cap
W\otimes k[t]$, $g_0 \in y^i I_{(E_i,f,x)}$, we shall be done.

Consider the (infinitely near) base points of the ideal $I_t$ for each
$t$ \cite[p. 254]{Cas00}. By hypothesis $\ell_E(1)\ge s$ and
$\hat\ell_E(0)\ge s$, so looking at the Newton polygon of elements in
$I_t$ we see that there 
are at least $s$ double base points on the germ $f_t=0$; on the other
hand $\ord(f_t|_D)\ge s$ so these base points do not depend on $t$ and
lie on the exceptional divisor. Let $\tilde \pi:S \rightarrow S_O$ be
the blowing up of the $s$ base points, and let $P\in S$ be the point
where the last 
exceptional divisor meets the strict transform $\tilde D$ of $D$.

To simplify matters and to be able to use the results
of section \ref{sec:formal}, we pass to the completion, as we may.
So let us denote $\O= \hat \O_{P,S}$ and let $\hat x, \hat y\in \O$ be local
equations of the last exceptional divisor and of $\tilde D$
respectively. We require in addition that $\tilde \pi^* f= \hat x^s(\hat x
+\hat y)\cdot u$, where $u$ is a unit (this is not restrictive, since
all ordinary singularities of multiplicicty three are analytically
equivalent). Then  $\hat x+\hat y+t$ differs from
$\tilde f_t$ by a unit, where $\tilde f_t=(\tilde \pi^* f_t)/\hat x^s$.
Let $\hat I_t$ be the completion of the virtual transform of $I_t$ 
with multiplicity two at the $s$ base points; i.e., 
$\hat I_t=\tilde I_t \otimes_{\O_{P,S}\otimes k[t]} (\O\otimes
k[[t]])$ 
where
$$
\tilde I_t=\frac{\tilde \pi^* (I_t)}{\hat x^{2s}} 
\subset\O_{P,S}\otimes k[t].
$$ 
With these notations, it is not hard to see that 
$$\hat I_t=I_{(\hat E,\hat x+\hat y+t,\hat x)},
$$
where $\hat E$ 
is obtained from $E$ 
by shortening the stair lengths by $s$, i.e., 
$\hat \ell_{\hat E}(i)=\hat \ell_{E}(j)-s$, $j=0,1.$

Denote by $V\subset \O_{P,S}\subset \O$ the set of virtual
transforms of elements of $W$ with (virtual) multiplicity 
at least 2 at the $s$ blown up points, i.e., 
$$V=\frac{\pi^*(\pi_*((\hat x^s \hat y)^m \hat x^{2s}))}{(\hat x^s
  \hat y)^m \hat x^{2s}} \otimes_{\O_{P,S}} \O,$$ 
where $\pi=\pi_O \circ \tilde \pi$ is the composition of the blow ups.

Remark that, by the proximity equality \cite[theorem 3.5.3]{Cas00}, every
$g\in V$ not multiple of $\hat y$ is the virtual transform at $P$ of a germ
at $O$ which has multiplicity exactly $m=2s$ and therefore does not
vanish at $P$ (i.e., $g \not\in (x,y)$). So if $T\subset \O/(y)$ is a
proper ideal of $\O/(y)$ then the canonical map
$$ 
\frac{V}{V \cap (y)} \overset{\varphi_p}{\longrightarrow}
\frac{R/(y)}{T}$$
is injective.

Now let $g_t \in I_t \cap (W\otimes k[t])$, and let $\hat g_t=\tilde
\pi^* (g_t)/\hat x^{2s}\in \hat I_t\cap(V\otimes k[[t]])$ be its
virtual transform. We have that $\hat \ell_{\hat E}(0)\ge 2$, so we may apply
proposition \ref{flatstairs} with $p=\ell_E(0)-2s$, which gives
$\tr{\hat I_t}{\hat y}{p}{}=1$ and then proposition \ref{formal} shows
that $\hat g_0 \in \hat y I_{(\sigma(\hat E,p),\hat f,\hat x)}$. So $g_0$
is a multiple of $y$: $g_0=yh$. Then $\hat g_0=\tilde \pi^*(g_0)/\hat
x^{2s}=\hat y \hat x^s \tilde \pi^*(h)/\hat x^{2s}$, with $\tilde
\pi^*(h)/\hat x^{s} \in I_{(\sigma(\hat E,p),\hat f,\hat x)}$. In
other words, $h\in \tilde \pi_* (\hat x^sI_{(\sigma(\hat E,p),\hat f,\hat
  x)})= I_{(E_1,f,x)}$, so in the case
$\ell_E(1)\le 2s$ we are done. 

Assume now that $\ell_E(1)\ge 2s+1$. In this case all elements $h\in
I_{(E_1,f,x)}$ as above have multiplicity at least $m+2$ along $D$. If
$y$ does not divide $h$, then $yh\in W$ tells us that $h$ is the
strict transform at $O'$ of a germ at $O$ which has multiplicity
exactly $m+1$, which again contradicts the proximity equality. So $y$
must divide $h$, and 
therefore $g_0 = y^2h'$ with $h'\in (I_{(E_1,f,x)}:y)=I_{(E_2,f,x)}$
and we are done.
\end{proof}

Theorem \ref{espbloc} is the main result on specialization inside the
Hilbert scheme that we shall use. As said above, our interest in the
schemes parameterized by the $H_{m,E,s}$ comes from the fact that they
lie in the border of the subscheme of $\Hilb \P^2$ parameterizing
2-curvilinear schemes. We now proceed to show this.

\begin{Lem}
  \label{lem:2curvmonom}
Let $Z$ be a zero-dimensional scheme supported at a single point
$O\subset \P^2$ and contained in a double curve $2C$, with $C$ smooth
at $O$. Let $N=\len Z$, $\ell=\len Z \cap C$, and let $E$ be
the staircase of height two with $\ell(0)=\ell$,
$\ell(1)=N-\ell$. Let $y=0$, $y \in \O_{O,\P^2}$ be a local equation
for $C$, and let $x \in \O_{O,\P^2}$ be transverse, so that $(x,y)$ is
the maximal ideal of $\O_{O,\P^2}$. Then there exists a flat family of
zero-dimensional schemes $Z_t \subset \P^2 \times \mathbb{A}^1$ such that
$Z_1=Z$, $Z_t$ is isomorphic to $Z$ for $t\ne 0$, and $Z_0$ is defined
by the ideal $I_{(E,y,x)}$.
\end{Lem}
\begin{proof}
Let $I\subset \O_{O,\P^2}\subset k[[x,y]]$ be the ideal defining
$Z$. As $I$ is $(x,y)$-primary, we may safely pass to the completion
$\hat \O_{O,\P^2}\cong k[[x,y]]$. It is immediate that $I_{(E,y,x)}$
is the initial ideal of $I$ with respect to the \emph{negative
  lexicographical ordering} with $1>x>y$ \cite[example
1.2.8]{GP02}. The desired family is then given by flat deformation to
the initial ideal (see, for instance, \cite[theorem 7.5.1]{GP02}).
\end{proof}

\begin{Lem}
  \label{lem:2curv1pt}
Let $Z$ be a zero-dimensional scheme contained in a double curve $2C$,
with $C$ smooth at $Z_{\red}$. Let $N=\len Z$, $\ell=\len Z \cap C$,
and let $E$ be the staircase of height two with $\ell(0)=\ell$,
$\ell(1)=N-\ell$. Let $O\in C$ be an arbitrary point, let $y=0$,
$y \in \O_{O,\P^2}$ be a local equation 
for $C$, and let $x \in \O_{O,\P^2}$ be transverse, so that $(x,y)$ is
the maximal ideal of $\O_{O,\P^2}$. Then there exists a flat family of
zero-dimensional schemes $Z_t \subset \P^2 \times \mathbb{A}^1$ such that
$Z_1=Z$, $Z_t$ is isomorphic to $Z$ for $t\ne 0$, and $Z_0$ is defined
by the ideal $I_{(E,y,x)}$.
\end{Lem}
\begin{proof}
Use \ref{lem:2curvmonom} and Hirschowitz's ``collision de front'' \cite{Hir85}.
\end{proof}

\begin{Cor}
    \label{lem:2curvesp}
Let $Z_0$ be a 2-curvilinear zero-dimensional scheme, with invariants
$N$ and $\ell$, and let 
$H(Z_0)\subset \Hilb \P^2$ be the set of all zero-dimensional subschemes
of the plane isomorphic to $Z_0$. 
For every positive integer $k$, let $E_k$ be the staircase of height
two and $\ell_{E_k}(0)=\ell-k(k+1)$, $\ell_{E_k}(1)=N-\ell-k^2$.
Then $H(Z_0)\subset \Hilb \P^2$ is constructible in the Zariski topology, and
for every $k\ge 1$ satisfying $2\ell -N > 2k$, $\ell \ge k(k+1)$, $N-\ell \ge
k^2$, one has $\overline{H(Z_0)}\supset H_{2k,E_k,k}$.
\end{Cor}
\begin{proof}
That $H(Z_0)$ is constructible is a general fact that
does not use 2-curvilinearity of $Z_0$. For the claimed incidences,
observe first that 
\ref{lem:2curv1pt} immediately gives 
$\overline{H(Z_0)} \supset H_{2,E_1,1}$, so for $k=1$ 
  we are done. Now proceed by recurrence on $k$,
  observing that for $k>1$, the hypotheses imply $\hat
  \ell_{E_{k-1}}(0)\ge k+1$ and $\ell_{E_{k-1}}(1)\ge 2k-1$, so
  theorem \ref{espbloc} tells us that $H_{2k,E_k,k}\subset
  \overline{H_{2k-2,E_{k-1},{k-1}}}$. 
\end{proof}

\section{Hilbert function of 2-curvilinear schemes}
\label{sec:2curv}

In order to prove our main theorem it is now enough to identify the
cases in which the sequence of specializations of the previous
section has led us to a maximal rank scheme type. To begin with, let
us recall a known class of maximal rank schemes:

\begin{Lem}
\label{wins}
  Let $E$ be a staircase of height two and $m$ a positive integer
  satisfying $\ell_E(0)> m$ and $2 \ell_E(1)\le m$. Then for every
  $s$ such that $m \ge \min \{e_1+se_2|(e_1,e_2)\in E\}$, general
  elements of $H_{m,E,s}$ have maximal rank.
\end{Lem}
\begin{proof}
  It is not hard to see that the schemes parameterized by $H_{m,E,s}$
  are cluster schemes (they are defined by complete ideals). Then it is
  enough to compute their cluster of base points to see that the claim
  is equivalent to lemma 4.4 of \cite{Roe01a}.
\end{proof}

\begin{proof}[Proof of theorem \ref{mainthm}]
  Note first that if $N-\ell \le 1$ then either $N=\ell$ and elements
  of $H_{m,E,s}$ are curvilinear, in which case the result is well
  known, or due to lemma \ref{lem:2curv1pt} it is enough to
  prove the maximal rank for general elements of $H_{2,E,\ell-2}$, with $E$ of
  height 1, which is again well known. See \cite[lemma 4.3]{Roe01a} for a
  proof that covers both cases over a field of arbitrary
  characteristic. Other proofs for the curvilinear case can be found in
  the literature; two elegant options over $\C$ are \cite{CM98}, which
  works in arbitrary dimension, or the use of Brian\c{c}on's 
  specializations of \cite{Bri77}.
 
  So assume that $N-\ell \ge 2$ and let $k$ be the maximal integer
  such that $N-\ell\ge k^2$. The hypothesis of the theorem tells us
  that $\ell \ge (N-\ell)+1+3\sqrt{N-\ell}\ge k^2+k$ and $2\ell -N
  \ge 1+3\sqrt{N-\ell}\ge 2k$, so we may apply corollary
  \ref{lem:2curvesp} and it will be enough to prove that general
  members of $H_{2k,E_k,k}$, with $E_k$ as in \ref{lem:2curvesp}, have
  maximal rank. We
  distinguish two cases. Assume first that $N-\ell \le k(k+1)$. Then
  it follows that $\ell_{E_k}(1)\le k$ and $\ell_{E_k}(0)\ge
  (N-\ell)-k^2+1+3\sqrt{N-\ell}-k\ge 2k+1$, and lemma \ref{wins}
  finishes the proof. 

Assume now that $N-\ell \ge k(k+1)+1$. Then
  $\ell_{E_k}(1) \ge k+1$ and $\hat \ell_{E_k}(0)\ge k+3$ so we may
  apply theorem \ref{espbloc} and obtain
  $\overline{H_{2k,E_k,k}}\supset H_{2k+1,E_k',k+1}$, where $E_k'$ has
  height two and $\ell_{E_k'}(0) =\ell -(k+1)^2-1$, $\ell_{E_k'}(1)
  =N-\ell- k(k+1)$. So it is enough to prove that general members of
  $H_{2k+1,E_k',k+1}$ have maximal rank. But $\ell_{E_k'}(0) =\ell
  -(k+1)^2-1\ge (N-\ell)-(k+1)^2+3\sqrt{N-\ell}>2k$ and $\ell_{E_k'}(1)
  =N-\ell- k(k+1)< k+1$ (because the choice of $k$ gives $N-\ell<
  (k+1)^2$) and again we finish using lemma \ref{wins}.
\end{proof}

\begin{proof}[Proof of theorem \ref{secondthm}]
 We want to apply theorem \ref{mainthm}. Clearly a scheme $Z_{(K,\m)}$
 is 2-curvilinear (as explained in the introduction) but we need some
 bounds on the length $N$ and maximal contact $\ell$ to hold.

   If $(D,\m)$ is the diagram of an $A_{2k-1}$ singularity, $k\ge 1$, and $K\in
  \Cl(D)$ has $p_1(K)=O\in \P^2$, let $(x,y)\in \O_{O,\P^2}$ be a
  system of parameters such that $y=0$ is the
  equation of a smooth germ of curve going through all the points of
  $K$. Then the ideal of $Z_{(K,\m)}$ is
  $(y^2,yx^k,x^{2k})$, and its invariants are $N=3k$,
  $\ell=2k$. Similarly, if $(D,\m)$ is the diagram of an 
  $A_{2k-2}$ singularity, $k\ge 2$, $K\in 
  \Cl(D)$ has $p_1(K)=O\in \P^2$, and $(x,y)\in \O_{O,\P^2}$ are a
  system of parameters such that $y=0$ is the
  equation of a smooth germ of curve going through all the points of
  $K$ but the last, then the ideal of $Z_{(K,\m)}$ is
  $(y^2,yx^k,x^{2k-1})$, and its invariants are $N=3k-1$,
  $\ell=2k-1$. All in all, the invariants of a union of singularity
  schemes always satisfy $5\ell \ge 3N$.

  Therefore a scheme as in the claim, being the union of such schemes,
  still satisfies $\ell \ge (3/5)N$. So in order to prove that $N
  \le 2 \ell -1 -3 \sqrt{N-\ell}$ it would be enough that $N \le
  (6/5)N-1-3\sqrt{(2/5)N}$. This inequality is always satisfied if
  $N\ge 100$ and then the claim follows from theorem
  \ref{mainthm}. Now assume that $N<100$; we need to prove the
  independence of the conditions imposed by a general $Z_{(K,\m)}$ to
  curves of degree 13 (and therefore of higher degree as
  well). Consider the scheme $X$ union of $Z_{(K,\m)}$ and $100-N$
  reduced points in general position. $X$ is still a 2-curvilinear
  scheme, and has 
  invariants $N'=100$, $\ell'=\ell+100-N\ge (3/5)/N'$ so it is of
  maximal rank as before. The linear system of all curves of degree 13
  has dimension $104>100$, so $X$ does impose independent conditions to
  curves of degree 13, and hence $Z_{(K,\m)}\subset X$ does too.
\end{proof}

Finally, let us deal with the low degree cases. First of all, there
are a number of (well known) exceptions.

\begin{Rem}
\label{exceptions}
  The following schemes are \emph{not} of maximal rank in the
  mentioned degrees:
  \begin{enumerate}
  \item Two ordinary double points in general position, in degree 2.
  \item Two ordinary cusp schemes in general position, in degree 3.
  \item A union of nodes or tacnodes $A_{2k_i-1}$ in general position
  with $\sum k_i=5$, in degree 4 (a particular case of which is 5
  double points mentioned in the introduction).
  \end{enumerate}
\end{Rem}

\begin{Teo}
\label{lastthm}
 Let $(D,\m)$ be a union of weighted diagrams of types $A_k$ (shown in
figure \ref{fig:Ak}), not among the exceptions \ref{exceptions}. 
Then for $K$ general in $\Cl(D)$, $Z_{(K,\m)}$
has maximal rank. 
\end{Teo}
\begin{proof}
Because of theorem \ref{secondthm}, we only need to show that the
given schemes have maximal rank in degrees less than 13. Let $N$ and
$\ell$ be the invariants associated to schemes $Z_{(K,\m)}$ for $K\in
\Cl(D)$. We can assume that $5\ell \ge 3N$. Let $d$ be the maximal
integer such that 
$d(d+1)/2<N$, and let $N'=d(d+1)/2$, $N''=(d+1)(d+2)/2$,
$\ell'=\min\{\ell,d(d+1)/2\}$, $\ell''=\ell+d(d+1)/2-N$ then it is
very easy to see that a scheme $Z \in H_{2,E(N,\ell),1}$ contains a 
$Z' \in \Hilb^{N'} \P^2$ and is contained in
a $Z'' \in \Hilb^{N''} \P^2$, and both are unions of singularity
schemes of multiplicity two and (possibly) simple points (in
particular, they are still 2-curvilinear and their invariants satisfy 
$5\ell \ge 3N$). Moreover, if $Z$ is in general position then we may
assume that $Z'$ and $Z''$ are in general position too. Thus,
reasoning as in \cite{Hir85} or \cite{Roe01a}, it is enough to prove
that every union of singularity schemes of multiplicity two and 
simple points, in general
position, of length $N=(d+1)(d+2)/2$ for $d \le 12$, has maximal rank.

The cases with $N \le 2\ell-1-3 \sqrt{N-\ell}$ have already been
solved. In particular for $N=(12+1)(12+2)/2$ all cases with $\ell\ge
55$ are done, and these are the only ones with $5\ell \ge 3N$, so we
can assume $d\le 11$. We consider the cases with $d\le 4$ to be well
known. Putting everything together the cases we are left with have the
following invariants:
$$
\begin{matrix}
  d& N & \ell \\
\hline 
11 &    78    & 47 \\
10 &    66 &    40 \\
9  &    55 &    33, 34 \\
8  &    45 &    27,28 \\
7  &    36 &    22, 23 \\
6  &    28 &    17,18 \\
5  &    21 &    13,14 \\
\end{matrix}
$$
The cases $N=55, \ell=33$ and $N=45, \ell=27$ can only be realized by
singularity schemes consisting of 11 and 9 ordinary cusps
respectively. Thus they were solved by Barkats \cite{Bar98}.

The case $N=78, \ell=47$ leads, after the sequence of specializations
used in the proof of theorem \ref{mainthm}, to schemes in
$H_{10,E,5}$, where $E$ is the height two staircase with
$\ell_E(0)=17$, $\ell_E(1)=6$. These are cluster schemes (defined by
complete ideals) which have maximal rank by \cite[proposition
4.5]{Roe01a}, so we are done. The same method works for the cases
$N=55, \ell=34$, $N=36, \ell=23$ and $N=21,\ell=14$.

The case $N=66, \ell=40$ can only be realized by singularity schemes
consisting of 10 ordinary cusps plus some other singularity
schemes. So a scheme $Z$ in this class can be written as $Y \cup Z'$
where $Y$ is a cusp scheme and $Z'$ is a 2-curvilinear singularity
scheme with invariants $N'=61, \ell'=37$, both in general
position. After the sequence of specializations used in the proof of
theorem \ref{mainthm}, $Z'$ degenerates into a scheme $Z'_0$ in
$H_{8,E,4}$, where 
$E$ is the height two staircase with $\ell_E(0)=17$, $\ell_E(1)=8$, so
$Z$ degenerates into $Y \cup Z'_0$. Specialize now the position of $Y$
so that its tangent line meets $Z'_0$, which forces curves of degree
10 containing $Y \cup Z'_0$ to contain the line as well. It is not hard
to check that in fact these curves consist of the line tangent to $Y$
counted twice plus curves of degree 8 through a residual scheme of
maximal rank and degree 45 ($Z'_0$ is in fact a cluster scheme so
residuals are easily computed) so there are no such curves. By
semicontinuity then there are no curves of degree 10 containing $Z$
either, and we are done. The same argument solves cases $N=45,
\ell=28$ and $N=36, \ell=22$ (which can only be realized by schemes
one of whose components is an ordinary cusp scheme). Not all cases with $N=28,
\ell=18$ have one ordinary cusp as a component, but those which do
have one are also solved.

For the remaining cases we only state the specializations that lead to
the solution, leaving to the reader the actual computations.

If $N=28,\ell=18$ and no component is an ordinary cusp, let $3\le k\le 5$
be minimal such that one component is the scheme of a singularity of
type $A_{2k-2}$ (analytically equivalent to $y^2-x^{k+1}$). If $k=3$,
specialize this component to a scheme in $H_{3,E_1,2}$, where $E_1$
has height 1 and length 2, and the rest of the components (as in
theorem \ref{mainthm})  to a scheme in $H_{4,E_2,2}$, where $E_1$
has height 2 and stairs of lengths 13 and 6. The line joining
the two components is a fixed part of the system under consideration
and allows to conclude. If $k=4$,
specialize the rest of the components (as in
theorem \ref{mainthm})  to a scheme in $H_{4,E,2}$, where $E$
has height 2 and lengths 5 and 2. The line joining
the two components is a fixed part of the system under consideration
and allows to conclude. If $k=5$, then there are exactly two
components both of type $A_6$. Specialize one of them (as in
theorem \ref{mainthm})  to a scheme in $H_{4,E,2}$, where $E$
has height 2 and lengths 3 and 2, and specialize further so that it is
supported on the unique conic of maximal contact with the other
component. This conic and the line tangent to the scheme of
$H_{4,E,2}$ are fixed parts of the system under consideration
and allow to conclude. 

If $N=28,\ell=17$, there must be 4 ordinary cusp schemes
involved. Specialize two of them to be supported at the tangent line
of a third; this line is a fixed component and allows to conclude.

Finally, if $N=21, \ell=13$ there must be at least one ordinary
cusp. If there are two, then the rest of the components can be
specialized to a scheme of a singularity of type $A_6$ with contact 6
with a line, which allows to conclude. If there is only one, then
one other component must be the scheme of a singularity of type $A_4$;
specialize so that the cusp is supported at the tangent line to the $A_4$.
\end{proof}

\bibliographystyle{amsplain}
\bibliography{Biblio}
\end{document}